\documentstyle{article}

\newtheorem{df}{ \sc Definition}[section]

\newtheorem{ex}[df]{ \it Example}
\newtheorem{pr}[df]{ \sc Proposition}

\newtheorem{tm}[df]{ \sc Theorem}
\newtheorem{cor}[df]{ \sc Corollary}
\newtheorem{re}[df]{ \it Remark}

\newtheorem{lem}[df]{ \sc Lemma}

\newtheorem{eq}[df]{\rm}

\def\hat{\widehat}
\def\bar{\overline}
\def\log{{\rm log}}
\def\setminus{\smash-\smash}
\def\p{\phantom{\mpd{}^\mpd{}_\mpd{}}}
\def\mpr#1{\;\smash{\mathop{\hbox to 20pt{\rightarrowfill}}\limits^{#1}}\;}
\def\mpl#1{\;\smash{\mathop{\hbox to 20pt{\leftarrowfill}}\limits^{#1}}\;}
\def\epi#1{\;\smash{\mathop{\hbox to 20pt{\rightarrowfill}\hskip
-15pt\rightarrow}\limits^{#1\,}}\;}
\def\mon#1{\;\smash{\mathop{\hbox to 
20pt{$\lhook\joinrel\relbar\joinrel\rightarrow$}}\limits^{#1\,}}\;}
\def\mono{\lhook\joinrel\relbar\joinrel\rightarrow}
\def\mpd#1{\big\downarrow\rlap{$\vcenter{\hbox{$\scriptstyle#1$}}$}}

\def\mpu#1{\big\uparrow\rlap{$\vcenter{\hbox{$\scriptstyle#1$}}$}}
\def\O{{\cal O}}

\font\Ch=msbm10

\def\Nt{\hbox{\Ch N}}

\def\Ct{\hbox{\Ch C}}
\def\Pt{\hbox{\Ch P}}

\def\s{{\rm sing}}
\def\o{{\rm reg}}
\def\Proof{\noindent\hskip4pt{\it Proof}}
\def\qed{\phantom a
\hfill$\vbox{\hrule\hbox{\vrule\vbox{\phantom{P}}\vrule}\hrule}$\vskip6pt}
\def\r{{Res(\omega)}}

\begin{document}
\date{February, 2005}
\author{Andrzej Weber\footnote{
Supported by KBN 1 P03A 005 26 grant.
I thank Institute of Mathematics, Polish Academy of Science for
hospitality.} }
\title{\bf Residue forms on singular hypersurfaces}
\maketitle
{\footnotesize
\noindent\hfil {Math. Sub. Class.: 14 F10, 14F43, 14C30}\vskip 10pt

}
\vskip 6pt
\section{Introduction}

The purpose of this paper is to point out a relation between the
canonical sheaf and the intersection complex of a singular algebraic
variety. We focus on the hypersurface case.
Let $M$ be a complex manifold, $X\subset M$ a
singular hypersurface. We study residues of top-dimensional
meromorphic forms with poles along $X$.
Applying resolution of singularities sometimes we are able to
construct residue
classes either in $L^2$-cohomology of $X$ or in the intersection
cohomology. The conditions allowing to construct these classes
coincide. They can be formulated in terms of the weight filtration.
Finally, provided that these conditions hold, we 
construct in a canonical way a lift of the residue class to cohomology of $X$.
\vskip 6pt

Let the manifold $M$ be of dimension $n+1$.
If the hypersurface $X$ is smooth we have an exact sequence of
sheaves on $M$:
$$0\mpr{}\Omega^{n+1}_M\mono \Omega^{n+1}_M(X)\epi{Res}
i_*\Omega^n_X\mpr{}0\,.$$
Here $\Omega^{n+1}_M$ stands for the sheaf of holomorphic differential
forms of the top degree on $M$ and $\Omega^{n+1}_M(X)$ is the sheaf of
meromorphic forms with logarithmic poles along $X$, i.e.~with the
poles at most of the first order.
The map $i:X\mono M$ is the inclusion. The morphism $Res$ is the residue map
sending $\omega={ds\over s}\wedge\eta$ to $\eta_{|X}$ if $s$ is a
local equation of $X$. The residues of forms with logarithmic poles
along a smooth hypersurface were studied by Leray (\cite{Le}) for
forms of any degree. Later such forms and their residues were applied
by Deligne (\cite{De}, \cite{GS})
to construct the mixed Hodge structure for the cohomology of open
smooth algebraic varieties.

We will allow $X$ to have singularities. As in the smooth case the
residue form is well defined differential form on the nonsingular
part of $X$. In general it may be highly singular at the singular
points of $X$. We will ask the following questions:
\begin{itemize}
\item Suppose $M$ is equipped with a hermitian metric. Is the norm of
$\r $ square integrable? We note that this condition does not depend on
the metric.
\item Does the residue form $\r $ define a class in the
intersection cohomology $IH^n(X)$?
\end{itemize}
We recall that by Poincar\'e duality residue defines a class in
homology $H_n(X)$ (precisely Borel-Moore homology, i.e.~homology with
closed supports), see \S\ref{rh}. The possibility to lift the residue class to
intersection cohomology means that $\r $ has mild
singularities. The intersection cohomology $IH^*(X)$, defined in
\cite{GM}, is a certain cohomology group attached to a singular
variety. Poincar\'e duality map $[X]\cap:H^*(X)\rightarrow
H_{\dim_{\bf R}(X)-*}(X)$
factors through $IH^*(X)$. Conjectu\-rally\footnote
{The proof in \cite{Oh} seems to be incomplete.}
 intersection cohomology is
isomorphic to $L^2$-cohomology. It was known from the very
beginning
of the theory, that conjecture is true if $X$ has conical
singularities (\cite{Ch}, \cite{CGM}). 

We study a resolution of singularities
$$\mu:\widetilde M\mpr{} M\,, \qquad \mu^{-1}(X)=\widetilde X\cup E\,,$$
where $\widetilde X$ is the proper transform of $X$ and $E$ is the
exceptional divisor. The pull back $\mu^*\omega$ is a meromorphic
form on $\widetilde M$. It can happen that it has no poles along the
exceptional divisors. Then we say that $\omega$ has {\it canonical}
singularities along $X$. By the definition $\omega$ has canonical
singularities if and only if $\omega\in adj_X\cdot\Omega^{n+1}_M(X)$,
where $adj_X\subset \O_M$ is 
the
adjoint ideal of {\rm \cite{EL}}.
The set of forms with canonical
singularities can be characterized as follows:

\begin{tm} \label{teo} The following conditions are equivalent
\begin{itemize}
\item $\omega$ has canonical singularities along $X$;
\item the residue form $\r\in\Omega^n_{X_\o}$ extends to a
holomorphic form on any resolution $\nu:\bar X\rightarrow X$;
\item the norm of $\r $ is square-integrable for any
hermitian metric on $X_\o$.
\end{itemize}
\end{tm}

Later the statement of \ref{teo} is divided into Proposition
\ref{p33}, Theorem \ref{hr}, Corollary \ref{c51} and Theorem 6.1.
Although our constructions use resolution of singularities we are
primarily interested in the geometry of the singular space $X$ itself.
The resulting objects do not depend on the choice of resolution.

Our description of forms with canonical singularities agrees with
certain results concerning intersection cohomology. We stress
that on
the level of forms we obtain a lift of residue to $L^2$-cohomology for
free. On the other hand, using cohomological methods one constructs a
lift of the residue class to intersection cohomology. This time the
lift is obtained essentially applying the decomposition theorem of
\cite{BBD}. This lift is not unique. It is worth-while to confront
these two approaches. The crucial notion in the cohomological
approach is the weight filtration. We will sketch this construction
below: Suppose
that $M$ is complete. Then $H^{k+1}(M\setminus X)$ is equipped with
the weight filtration, all terms are of the weight $\geq k+1$.
The homology $H_{2n-k}(X)$ is also 
equipped with a mixed Hodge structure. It is of the
weight $\geq k-2n$. The homological residue map preserves the weight
filtration: 
$$res:H^{k+1}(M\setminus X)\mpr{}H_{2n-k}(X)(-n-1).$$
Here $(i)$ denotes the $i$-fold Tate twist; now $H_{2n-k}(X)(-n-1)$ is of
the weight $\geq k+2$.
Intersection cohomology $IH^k(X)$ maps to $H_{2n-k}(X)(-n)$. Since it is
pure, the image is contained in $W_k(H_{2n-k}(X)(-n))$.
We will show that:

\begin{tm} \label{teo1}If $c\in H^{k+1}(M\setminus X)$ is of weight $\leq k+2$
then $res(c)$ lifts to intersection cohomology. In another words, we
have a factorization of the residue map
$$\def\mpr#1{\;\smash{\mathop{\hbox to 40pt{\rightarrowfill}}\limits^{#1}}\;}
\matrix{W_{k+2}H^{k+1}(M\setminus X)&\mpr{res}&W_{k+2}(H_{2n-k}(X)(-n-1))&.\cr
\hfill\searrow&&\nearrow\hfill\cr
&IH^k(X)(-1)&\cr}$$
\end{tm}
 In fact for an arbitrary complete algebraic variety the image of
intersection cohomology coincides with the lowest term of the
weight filtration in homology, see \cite{We4}.

We note that if $\omega$ has canonical singularities along $X$ then
its cohomology class is of weight $\leq n+2$. By \ref{teo}
$\r $ defines a class in $L^2$-cohomology. Also, by
\ref{teo1} the residue of $[\omega]$ can be lifted to intersection
cohomology. To completely clear up this situation we construct in
\S\ref{rich} a canonical lift of the residue class not only to
intersection cohomology, but even to cohomology of $X$. 

An attempt to relate holomorphic differential forms to
intersection cohomology was proposed by Koll\'ar (\cite{Ko1}, II
\S4 ). It seems that his solution is not definite since he
applies the (noncanonical) decomposition theorem. The construction
proposed in \ref{ltc} is elementary and geometric. As a side
result of our consideration we obtain

\begin{tm}\label{wyz} Suppose an algebraic variety $X$ is complete of 
dimension
$n$. Let $\widetilde X$ be its resolution. Then $H^k(\widetilde
X;\Omega_{\widetilde X}^n)$ is a direct summand both in
$H^{n+k}(X)$ and $IH^{n+k}(X)$. \end{tm}

One can hope that a relation between holomorphic forms of lower
degrees with intersection cohomology will be explained as well.

Another approach to understand the relation between the residues
and intersection cohomology was presented by Vilonen \cite{Vi} in
the language of $\cal D$-modules. His method applies to isolated
complete intersection singularities.

Finally in \S\ref{oi}-\S\ref{qs} we briefly describe a relation
between the oscillating integrals of \cite{Ma} or \cite{Va} and
residue theory for isolated singularities. Namely, if the order of
a form at each singular point is greater than zero, then the
residue class can be lifted to intersection cohomology. Again,
this condition coincides with having canonical singularities.

The present paper is a continuation of \cite{We2}, where the case of
isolated singularities was described.
My approach here was partially motivated by a series of lectures
delivered by Tomasz Szemberg on the algebraic geometry seminar IMPANGA in
Polish Academy of Science.

{\def\contentsname{\normalsize\sc Contents:}
\def\bfseries{\rm}
\tableofcontents}

\section{Residues as differential forms}

Let $\omega$ be a closed form with a first order pole on $X$. Then
the residue form $\r $ can be defined at the regular points
of $X$. The case when $\omega$ is a holomorphic $(n+1,0)$--form is
the most important for us:
$$\omega=\frac gs dz_0\wedge\dots\wedge dz_n\,,$$
where the function $s$ describes $X$. The space of such
forms is denoted by $\Omega^{n+1}_M(X)$.
Then the residue form is a
holomorphic $(n,0)$--form:
$$\r \in \Omega^n_{X_\o}\,.$$
The symbol $\in$ by abuse of notation means $\r $ is a
section of the sheaf $\Omega^n_{X_\o}$.
The precise formula for the residue is the following: Set $s_i={\partial
s\over \partial z_i}$. We have
$$ds=\sum_{i=0}^{n}{s_i}dz_i\,.$$
At the points where ${s_0}\not=0$ we write
$$dz_0=\frac1{s_0}\left(ds-\sum_{i=1}^{n}{s_i}dz_i\right)$$ and
$$\omega={g\over {s\,s_0}}\left(ds-\sum_{i=1}^{n}{s_i}
dz_i\right)\wedge dz_1\wedge\dots\wedge dz_{n}=$$
$$={ds\over s}\wedge \frac g{s_0} dz_1\wedge\dots\wedge dz_{n}\,.$$
Thus
$\r = \left(\frac g{s_0} dz_1\wedge\dots\wedge
dz_{n}\right)_{|{X_\o}}\in \Omega^n_{X_\o}$.

To see how $\r$ behaves in a neighbourhood of the singularities
let us calculate its norm in the metric coming from the coordinate
system: $$|\r |_X=\left|{ds\over |ds|}\wedge\r\right|_M=
\left|{s\;\omega\over |ds|}\right|_M= {|g|\over |{\rm
grad}(s)|}\,.$$ We conclude that $\r $ has (in general) a pole at
singular points of $X$.

The forms that can appear as residue forms are exactly the regular
differential forms defined by Kunz for arbitrary varieties;
\cite{Ku}.

\section{Residues and resolution}
\label{roz}

We will analyze the residue form using resolution of singularities.
Let $\mu:\widetilde M\rightarrow M$ be a
log-resolution of $(M,X)$, i.e.~a birational
map , such that $\mu^{-1}X$ is a
smooth divisor with normal crossings and $\mu$ is an isomorphism when
restricted to
$\widetilde M\setminus \mu^{-1}X_{\rm sing}$. Let $\widetilde X$ be
the proper transform of $X$ and let $E$ be the exceptional divisor.
The pull-back of $\omega$ to $\widetilde M$ is a meromorphic form
with poles along $\widetilde X$ and $E$. According to the terminology
of \cite{Ko2} we define:

\begin{df} We say that $\omega$ has canonical singularities along $X$
if $\mu^*\omega$ has no pole along the exceptional divisor,
i.e. $\mu^*\omega\in\Omega^{n+1}_{\widetilde M}(\widetilde X)$.
\end{df}

We note that this notion does not depend on the resolution.
Our method of studying residue forms are appropriate to tackle this
class of singularities.
We begin with an easy observation:

\begin{pr} \label{p33} If $\omega$ has canonical singularities along $X$, then
for
any re\-so\-lution $\nu:\bar X\rightarrow X$ the pull-back of the residue
form $\nu^*\r $ is holomorphic on $\bar X$.\end{pr}

\begin{re}\rm We do not assume that $\nu$ extends to a resolution of the
pair $(M,X)$.\end{re}

\Proof. Let $\mu$ be a log-resolution of $(M,X)$. By the assumption
$\mu^*\omega\in\Omega^{n+1}_{\widetilde M}(\widetilde X)$. Therefore
$Res(\mu^*\omega)$ is a holomorphic form on $\widetilde X$. Hence
$\r\in\Omega^n_{X_\o}$ extends to a section of
$\mu_*\Omega^n_{\widetilde X}$. The later sheaf does not depend on
the resolution of $X$. Indeed, let $\hat X$ be a smooth variety
dominating both $\widetilde X$ and $\bar X$. Then $Res(\mu^*\omega)$
can be pulled back to $\hat X$ and pushed down to $\bar X$ (since
$f_*\Omega^n_{\hat X}=\Omega^n_{\bar X}$ if $f$ is birational). The
resulting form coincides with $\nu^*\r $ outside the
singularities. \qed

\section{Vanishing of hidden residues}

We have observed that if $\omega$ has canonical singularities, then
the residue form is smooth on each resolution. Let us assume the
converse: suppose $\r $ extends to a holomorphic form on
$\widetilde X$. The extension is determined only by the nonsingular part
of $X$. We will show, that all the other ''hidden'' residues along
exceptional divisors vanish.

\begin{tm}\label{hr} If $Res(\mu^*\omega)_{|\widetilde X\setminus E}$ has no
pole along $E\cap\widetilde X$ then $\omega$ has canonical singularities along
$X$.
\end{tm}

\Proof. Let $E=\bigcup_{i=1}^kE_i$ be a decomposition of $E$ into
irreducible components.
 Assume that $Res(\mu^*\omega)_{|E_i}$ is nontrivial for $1\leq
i\leq l$ for some $l\leq k$. Blowing up intersections $E_i\cap
\widetilde X$ we can assume that $E_i\cap
\widetilde X=\emptyset$ for $i\leq l$. Let $a_i$ be the order of the
pole of $\mu^*\omega$ along $E_i$. Define a quotient sheaf $\cal F$:
\begin{eq}\hfil\label{ss}
$0\rightarrow \Omega^{n+1}_{\widetilde M}\mono
\Omega^{n+1}_{\widetilde M}(\sum_{i=1}^la_iE_i)\epi{}{\cal
F}\rightarrow 0\,.$\end{eq}

\begin{lem} \label{ll}The direct image $\mu_*\cal F$ vanishes.\end{lem}

\Proof. We push forward the sequence \ref{ss} and we obtain again the
exact sequence, since $R^1\mu_*\Omega^{n+1}_{\widetilde M}=0$ e.g.~by
\cite{Ko1}. But now the sections of $$\mu_*\Omega^{n+1}_{\widetilde
M}(\sum_{i=1}^la_iE_i)$$ are forms which are holomorphic on
$M\setminus\mu(E)$. Therefore they are holomorphic and hence
$\mu_*{\cal F}=0$.\qed

\Proof { \it of \ref{hr} cont.} We tensor the sequence
\ref{ss} with $\O(\widetilde X)$. Since the support of $\cal F$ is
disjoint with $\widetilde X$ we obtain a short exact sequence:
$$0\rightarrow \Omega^{n+1}_{\widetilde M}(\widetilde
X)\rightarrow \Omega^{n+1}_{\widetilde M}(\widetilde X+
\sum_{i=1}^la_iE_i)\rightarrow {\cal F}\rightarrow 0\,.$$ We apply
$\mu_*$ and by the Lemma \ref{ll} we have an isomorphism
$$\mu_*\Omega^{n+1}_{\widetilde M}(\widetilde X)\simeq
\mu_*\Omega^{n+1}_{\widetilde M}(\widetilde
X+\sum_{i=1}^la_iE_i)\,.$$ The above equality means that $\omega$
cannot have a pole along exceptional divisors.
\qed

\section{Adjoint ideals}

The adjoint ideals were introduced in \cite{EL} for a hypersurface
$X\subset M$. The adjoint ideal $adj_X\subset \O_M$ is the ideal satisfying
$$\mu_*\Omega^{n+1}_{\widetilde M}(\widetilde X)=
adj_X\cdot\Omega^{n+1}_M(X)\,.$$
 
The ideal $adj_X$ consists of the functions $f$, for which $\mu^*(\frac fs
dz_1\wedge\dots\wedge dz_m)\in \Omega^{n+1}_{\widetilde
M}(\mu^*D)$ has no pole along the exceptional divisors, i.e.~it
belongs to $\Omega^{n+1}_{\widetilde
M}(\widetilde X)$. Here $s$, as before, is a function describing
$X$. In another words the forms $\omega\in
adj_X\cdot\Omega^{n+1}_{ M} (X)$ are exactly the forms with
canonical singularities along $X$. Moreover the sequence of sheaves
\begin{eq}\label{es}\hfil $0\rightarrow\Omega^{n+1}_M \rightarrow
adj_X\cdot\Omega^{n+1}_M (X)
\rightarrow\mu_*\Omega^n_{\widetilde X}\rightarrow 0 $\end{eq}
is exact (\cite{EL} 3.1). 
 (This follows from vanishing of
$R^1\mu_*\Omega^{n+1}_{\widetilde M}$.)
The adjoint ideal does not depend on the resolution.

\begin{cor} \label{c51} The residue form $\r\in\Omega^n_{X_\o}$ extends
to a section of $\mu_*\Omega^n_{\widetilde X}$ if and only if $\omega\in
adj_X\cdot\Omega^{n+1}_M (X)$.\end{cor}

\Proof. The implication $\Rightarrow$ follows from the Theorem
\ref{hr}. The converse follows from the exact sequence \ref{es}.
\qed

It turns out that every form has canonical singularities, i.e.
$adj_X=\O_M$ if and only if $X$ has rational singularities \cite{Ko2}, \S11.

\section{$L^2$-cohomology }

Let us assume that the tangent space of $M$ is equipped with a
hermitian metric. For example if $M$ is a projective variety, then
one has the restriction of the Fubini-Study metric from projective
space. The nonsingular part of the hypersurface $X$ also inherits
this metric. One considers the complex of differential forms which have
square-integrable pointwise norm (and the same holds for
differential). Its cohomology is an important invariant of the
singular variety called the
$L^2$-cohomology, \cite{CGM}. This is why we are led to
the question: when the norm of the residue form is
square-integrable? Moreover, for the forms of the type $(n,0)$ on
the $n$-dimensional manifold the condition of integrability does not
depend on the metric. This is because $\int_X|\eta|^2dvol(X)$ is
equal up to a constant to $\int_X\eta\wedge\bar\eta$.

\begin{tm} The residue form $\r $ has the square-integrable
norm if and only if $\omega$ has canonical singularities.\end{tm}

\Proof. Instead of asking about integrability on $X_\o$ we ask about
integrability on $\widetilde X$. Now, local computation shows that if
$\omega$ has a pole, then its norm is not square-integrable.
\qed

\begin{re}\rm Note that the class of the residue form does not
vanish provided that $\omega$ has a pole along $X$. This is
because
$\r $ can be paired with its conjugate
$\overline{\r }$ in cohomology.\end{re}

\begin{re}\rm The connection between integrability condition and
multiplicities were studied by Demailly, see e.g.~\cite{Dm}.
\end{re}

\begin{re}\rm For homogeneous singularities
(which are conical) integrals of the residue forms along conical cycles
converge provided that the cycle is allowable in the sense of
intersection homology and $|\r |\in L^2(X)$.
\end{re}

\section{Residues and homology}
\label{rh}

Suppose for a moment that $X\subset M$ is smooth.
Let $Tub_X$ be a tubular neighbourhood of $X$ in $M$.
We have a commutative diagram: 
$$\matrix{ H^*(M\setminus X)&\mpr{[M]\cap}&H^{BM}_{2n+2-*}(M,X)(-n-1)\cr\cr
\mpd{d}&&\mpd{\partial}\cr\cr 
H^{*+1}(M,M\setminus X) &\mpr{[M]\cap}&H^{BM}_{2n+1-*}(X)(-n-1)\cr\cr
\parallel&&\mpu{[X]\cap}\cr\cr
H^{*+1}(Tub_X,Tub_X \setminus X)&\mpl{\tau}&H^{*-1}(X)(-1)}$$
In the diagram $H^{BM}_*$ denotes Borel--Moore
homology, i.e. homology with closed supports.
All coefficients are in $\Ct$. The entries of the diagram are equipped
with the Hodge structure. 
The map $\tau$ is the Thom isomorphism, the remaining maps in
the bottom square are also isomorphisms by Poincar\'e duality for $X$
and $M$.
The residue map
$$res=\tau^{-1}\circ d: H^*(M\setminus X)\longrightarrow H^{*-1}(X)$$
is defined to be the composition of the differential with the inverse of the
Thom isomorphism. By \cite{Le} we have:
$$res([\omega])={1\over 2\pi i}[\r ]$$
for a closed form with the first order pole along $X$.
(We use small letter for the homology class $res(c)\in H_{2n+1-*}^{BM}(X)$
to distinguish it from the differential form $\r\in\Omega^*_{X}$.)

When $X$ is singular then there is no tubular
neighbourhood of $X$ nor Thom isomorphism, but we can still define
a homological residue
$$res: H^*(M\setminus X)\longrightarrow H_{2n+1-*}^{BM}(X)(-n-1)$$
$$res(c) = [M]\cap dc=\partial([M]\cap c)$$
If $X$ was nonsingular, then this definition would be equivalent to the
previous one since $\xi\mapsto [X]\cap\xi$ is Poincar\'e duality
isomorphism and
the diagram above commutes. 
\begin{re} \rm One should mention here the work of M.~Herrera \cite{He1} and
\cite{HeL}  who defined
a residue current for a meromorphic $k+1$--form. This current is
supported by the divisor of poles. For a closed form it defines a
homology class in $H^{BM}_{2n-k}(X)$.\end{re}

In general there is no hope to
lift the residue morphism to cohomology. For $M=\Ct^{n+1}$ the morphism
$res$ is the Alexander duality isomorphism and $[X]\cap$ may be
not onto. Instead we ask if the residue of an element lifts to
the intersection homology of $X$.
The intersection homology groups, defined by Goresky and MacPherson
in \cite{GM}, are the groups that 'lie between'
homology and cohomology; i.e. there is a factorization:
 
$$\def\mpr#1{\;\smash{\mathop{\hbox to 30pt{\rightarrowfill}}\limits^{#1}}\;}
\matrix{
H^{*}(X) & \mpr{[X]\cap} &H_{2n-*}^{BM}(X)(-n)\,.\cr
\hfill\searrow & &\nearrow\hfill \cr
& IH^*(X)&\cr}
$$
 In fact for complete $X$ the map $[X]\cap$ factors through
$$ H^k(X)/W_{k-1}H^k(X)\mon{\alpha} IH^k(X) \epi{\beta}
W_k(H_{2n-k}(X)(-n))\,.$$
The injectivity of $\alpha$ and surjectivity of $\beta$ is proved in
\cite{We4}. The composition $\beta\alpha$ does not have to be an
isomorphism. For example, if $X$ admits an algebraic cellular
decomposition then 
its cohomology is pure (i.e. $W_{k-1}H^k(X)=0$ and
$W_k(H_{2n-k}(X)(-n))=H_{2n-k}(X)(-n)\,$) but the Poincar\'e
duality map $[X]\cap$ does not have to be an isomorphism.
We will analyze the arguments of \cite{We4} for the particular situation of
a hypersurface.

\section{Hodge theory}\label{lc}
According to Deligne (\cite{De}, see also \cite{GS}) any algebraic
variety carries a mixed Hodge structure. Suppose the ambient
variety $M$ is complete. To construct the mixed Hodge structure on
$M\setminus X$ one finds a log-resolution of $(M,X)$, denoted by
$\mu:\widetilde M\rightarrow M$, (see \S\ref{roz}). Then one
defines $A_{\rm log}^*=A^*_{\widetilde
M}(\log\langle\mu^{-1}X\rangle)$, the complex of $C^\infty$ forms
with logarithmic poles along $\mu^{-1}X$. Its cohomology computes
$H^*(\widetilde M\setminus\mu^{-1}X)=H^*(M\setminus X)$. The
complex $A_{\log}$ is filtered by the weight filtration
$$0=W_{k-1}A^k_\log\subset W_{k}A^k_\log\subset
\dots\subset W_{2k}A^k_\log=A^k_\log\,,$$
which we describe below. Let $z_0,z_1,\dots z_n$ be local coordinates
in which the components of $\mu^{-1}X$ are given by the equations
$z_i=0$ for $i\leq m$. The space $W_{k+\ell}A^k_\log$ is spanned by the forms
$$\frac{dz_{i_1}}{z_{i_1}}\wedge\dots
\wedge\frac{dz_{i_\ell}}{z_{i_\ell}}\wedge\eta$$
where $i_j\leq m$ and $\eta\in A^{k-\ell}_{\widetilde M}$ is a smooth
form on $\widetilde M$.
The weight filtration in $A^*_\log$ induces a filtration in
cohomology. The quotients of subsequent terms $W_{k+\ell}H^k(M\setminus
X)/W_{k+\ell-1}H^k(M\setminus X)$ are equipped with pure Hodge structure
of weight $k+\ell$.

Our goal is to tell whether the residue of a differential
form or the residue of a cohomology class can be lifted to
intersection cohomology. The Hodge structure on intersection
cohomology hasn't been constructed yet in the setup of differential
forms. On the other hand, there are alternative constructions in
which intersection homology has weight filtration. If $X$ is a
complete variety, then $IH^*(X)$ is pure. This property is
fundamental either in \cite{BBD} or in Saito's theory, \cite{Sa}.

The homology of $X$ is also equipped with the mixed Hodge structure.
Since $X$ is complete 
$$W_{k-1}(H_{2n-k}(X)(-n))=0\,,\qquad W_{2k}(H_{2n-k}(X)(-n))=
H_{2n-k}(X)(-n)\,.$$
Due to purity of intersection cohomology
$$im(IH^k(X)\rightarrow H_{2n-k}(X))\subset W_k(H_{2n-k}(X)(-n))\,.$$
The residue map 
$$res:H^{k+1}(M\setminus X)\rightarrow H_{2n-k}(X)(-n-1)$$
preserves the weights.
In particular it vanishes on
$$W_{k+1}H^{k+1}(M\setminus X)=im(H^{k+1}(\widetilde M)\rightarrow
H^{k+1}(\widetilde M \setminus \mu^{-1}X))\,.$$
Suppose we have a class $c\in W_{k+2}H^{k+1}(M\setminus X)$. Then
$res(c)$ is of weight $k+2$ in $H_{2n-k}(X)(-n-1)$. It is reasonable to ask
if it comes from intersection cohomology.

\begin{tm} \label{lift} Suppose that $M$ is complete. Then the
residue of each class $c\in
W_{k+2}H^{k+1}(M\setminus X)$ can be lifted to intersection
cohomology.\end{tm}

\Proof. Let $\mu:\widetilde M\rightarrow M$ be a log-resolution of
$(M,X)$. We consider the residue $res(\mu^*c)\in
H_{2n-k}(\mu^{-1}X)(-n-1)$.

\begin{lem} The homology class $res(\mu^*c)$ is a lift of res(c) to
$H_{2n-k}(\mu^{-1}X)(-n-1)$, i.e.
$$\mu_*(res(\mu^*c))=res(c)\,.$$\end{lem}

\Proof.
$$\mu_*(res(\mu^*c))=\mu_*([\widetilde M]\cap
d\mu^*c)=(\mu_*[\widetilde M])\cap dc=res(c)\,.$$
\qed

\Proof { \it of \ref{lift} cont.} Now assume that $c$ has
weight $k+2$. Then $\mu^*c$ is represented by a form $\omega$
with logarithmic
poles of weight $k+2$. The residue of $\omega$ consists of forms
$Res_i(\omega)$ on each component $E_i\subset \mu^{-1}X$ (we set
$E_0=\widetilde X$).
These forms have no poles along the intersections of components. This
means that $res(\mu^*c)$ comes from
$\sum_\alpha[Res_i(\omega)]\in \bigoplus_i H^k(E_i)=
IH^k(\mu^{-1}X)$. By \cite{BBFGK} (see \cite{We3} for a short proof) we can 
close
the following diagram with a map $\theta$ of intersection cohomology
groups:
$$\matrix{\sum_i[Res_i(\omega)]&\in&IH^k(\mu^{-1}X)&\mpr{\iota}
&H_{2n-k}(\mu^{-1}X)&\ni&res(\mu^*c)\cr
&&\mpd{\theta}&\p&\mpd{\mu_*}\cr
&&IH^k(X)&\mpr{\iota}&H_{2n-k}(X)&\ni&res(c)\,.\cr}$$
Here $\iota$ is the natural transformation from intersection
cohomology to homology.
The class $\theta(\sum_i[Res_i(\omega)])$ is the desired lift
of $res(c)$.
\qed

\begin{re}\rm The completeness assumption can be removed in
\ref{lift} and it is clear that the orders of poles at infinity do not
matter. \end{re}

Note that if a meromorphic $n+1$--form $\omega$ has canonical
singularities along $X$ then 
$\mu^*\omega$ has no pole along the exceptional divisors.
Therefore it belongs to the logarithmic complex, it is closed and
$$\mu^*\omega\in W_{n+2}A^{n+1}_\log\,.$$
Conversely, a closed $n+1$--form which belongs to the top piece of the Hodge
filtration $F^{n+1}A^{n+1}_\log$ has to be meromorphic. 
We obtain a surjection 
$$\Omega^{n+1}_{\widetilde M}(\mu^{-1}X)\epi{}F^{n+1}H^{n+1}(M\setminus 
X)\,.$$
A meromorphic form has canonical
singularities if and only if it belongs to $W_{n+2}A^{n+1}_\log$
since by
\ref{hr} it has no poles along the exceptional divisors. Therefore we
have a surjective map
$$\Omega^{n+1}_{\widetilde M}(\widetilde
X)\epi{}F^{n+1}W_{n+2}H^{n+1}(M\setminus X)\,.$$ 
This way we have solved positively the problem of lifting to $IH^n(X)$
the residue classes of forms which have canonical singularities. 
Nevertheless it is possible to do much more.
We will find a lift to cohomology $H^n(X)$ in a canonical way.

\section{\label{rich}Residues in cohomology }

In this section we ignore the Tate twist.

Suppose a meromorphic $(n+1)$-form $\omega$ has canonical
singularities along $X$. We will show how to construct a lift of
the residue class $res(\omega)\in H_n(X)$ to $H^n(X)$. It is
enough to define an integral
$$\widehat{res}(\omega): H_n(X)\rightarrow
\Ct\,.$$ For the construction we need the following (probably well known)
fact.

\begin{pr} Let $X$ be a variety of pure dimension. Let
$TC^{alg}_*(X)\subset C_*(X)$ be the subcomplex of geometric chains which
are semialgebraic and satisfy the conditions
$$\dim(\xi \cap X_\s)<\dim \xi\,,$$
$$\dim(\partial\xi \cap X_\s)<\dim \partial\xi\,.$$
The inclusion of complexes induces an isomorphism of homology.\end{pr}

\begin{re}\rm To show that the support condition does not spoil the
homology one can proceed as in \cite{Ha} computing inductively local
cohomology.\end{re}

For a cycle $\xi\in TC^{alg}_n(X)$ let us define 
$$\langle \widehat{res}(\omega),\xi\rangle = {1\over 2\pi i}
\int_{\mu^*\xi}Res(\mu^*\omega)\,,$$
where
$$\mu^*\xi={\rm closure}(\mu^{-1}(\xi \setminus X_\s))$$
is the strict transform of the cycle $\xi$. Note that $\mu^*\xi$
is a semialgebraic chain, which does not have to be a cycle.
Alternatively, we may define $\langle
\widehat{res}(\omega),\xi\rangle ={1\over 2\pi i }\int_{\xi}Res(\omega)$
and say that the integral always converges for $\xi\in
TC^{alg}_n(X)$. We have to prove, that our definition does not
depend on the choice of a cycle. Suppose that $\xi'$ is another
cycle, such that $\xi-\xi'=\partial \eta$. Again we assume that
both $\xi'$ and $\eta$ belong to $TC^{alg}_*(X)$. Set
$$\Delta=\mu^*\xi-\mu^*\xi'-\partial\mu^*\eta\,.$$
The residue form $Res(\mu^*\omega)$ is closed, therefore by Stokes
theorem
$$\matrix{\langle \widehat{res}(\omega),\xi\rangle - \langle
\widehat{res}(\omega),\xi'\rangle & = &
{1\over 2\pi i}\left(\int_{\mu^*\xi}Res(\mu^*\omega) -
\int_{\mu^*\xi'}Res(\mu^*\omega) \right)\cr
& = &
{1\over 2\pi i}\int_\Delta Res(\mu^*\omega)\,.\hfill}$$
The chain $\Delta$ is contained in the exceptional locus of
$\mu_{|\widetilde X}$, which is of dimension $n-1$. The form 
$Res(\mu^*\omega)$ is of
type $(n,0)$, therefore it vanishes on $\Delta$.
This way we have defined $\widehat{res}(\omega)\in (H_n(X))^*=H^n(X)$.

We have to show that $\widehat{res}(\omega)$ is a lift of
$res(\omega)\in H_n(X)$. In fact we will argue that it is a lift of
$res(\mu^*\omega)\in H_n(\mu^{-1}X)$.
By our assumption $res(\mu^*\omega)$ comes from $\bigoplus_i H^n(E_i)$.
By \ref{hr} the residues $Res_i(\mu^*\omega)$ vanish along
the exceptional divisors. It is enough to show that
$$\langle\widehat{res}(\omega),[\mu_*(\xi)]\rangle=
{1\over 2\pi i}\langle Res_0(\mu^*\omega),\xi\rangle=
{1\over 2\pi i}\int_\xi Res_0(\mu^*\omega)$$
for a cycle $\xi \in C_n(\widetilde X)$. We may assume that $\xi$ is
semialgebraic and $\dim(\xi\cap \mu^{-1}(X_\s)\leq n-1$. Then
$\mu^*\mu_*\xi=\xi$ and the formula follows from the definition of
$\widehat{res}(\omega)$.

We have proved
\begin{tm} If $\omega$ is a holomorphic form of the top degree, then there
exists a canonical lift of $res(\omega)$ to cohomology $H^n(X)$.\end{tm}

\begin{re}\rm \label{ltc} By the same procedure one can define a map
$$\iota:H^k(\widetilde X,\Omega^n_{\widetilde X})\rightarrow H^{n+k}(X)\,,$$
such that $\mu^*\circ\iota$ is the canonical map $H^k(\widetilde
X,\Omega^n_{\widetilde X}) \rightarrow H^{n+k}(\widetilde X)$. By
\cite{BBFGK} the map $\mu^*:H^*(X)\rightarrow H^*(\widetilde X)$
factors through $IH^*(X)$. On the level of derived category
$D(X)$ we have a chain of maps
$$R\mu_*\Omega^n_{\widetilde X}[-n] \simeq \mu_*A^{n,*}[-n]\rightarrow
\Ct_X\rightarrow IC_X \rightarrow R\mu_*\Ct_{\widetilde X}$$
factorizing the natural $R\mu_*\Omega^n_{\widetilde
X}[-n]\rightarrow R\mu_*\Ct_{\widetilde X}$. 
This proves Theorem \ref{wyz}.
Note that a map to
intersection cohomology or rather a dual one
$$IC_X\rightarrow DR\mu_*\Omega^n_{\widetilde X}[-n]\simeq R\mu_*{\cal 
O}_{\widetilde X}$$
was described by Koll\'ar in \cite{Ko1}, II 4.8. The decomposition
theorem of \cite{BBD} is applied. Our map is constructed
surprisingly easily and in a canonical way.

For complete $X$ we obtain a side result:
\end{re}

\begin{tm} Suppose an algebraic variety $X$ is complete of dimension
$n$. Let $\widetilde X$ be its resolution. Then $H^k(\widetilde
X;\Omega_{\widetilde X}^n)$ is a direct summand both in
$H^{n+k}(X)$ and $IH^{n+k}(X)$. The inclusion is adjoint to the
strict transform of cycles.\end{tm}

The statement for intersection cohomology also follows from \cite{Ko1}
II 4.9.

\begin{re}\rm In \cite{He2} there are studied residues of the meromorphic 
forms which
can be written as $\omega={ds\over s}\wedge\eta+\theta$.  For the
forms of top degree this condition is more restrictive then having canonical
singularities. 
For example if $n\geq 2$ and $X$ has isolated simple
singularities then all forms $\omega\in\Omega^{n+1}(X)$ have
canonical singularities 
(see \S11) but not necessarily can be written as above.
For the forms considered by Herrera the residue
$res(\omega)=\eta_{|X}$ is well defined as an element of a suitable
complex of forms on the singular variety $X$. 
The space $M$ is allowed to be singular. For nonsingular $M$ this
result is rather tautological.\end{re}

\section{Isolated singularities }

\label{oi} Residue forms for hypersurfaces with isolated
singularities are strongly related to oscillating integrals. The
first references for this theory are \cite{Ma} or \cite{Va}. In
\cite{AGVII}\S10-15 the reader can find a review, samples of proofs and
other precise references to original papers. A relation of
oscillating integrals with the theory of singularities of pairs is
explained in \cite{Ko2}, \S9.

Suppose $0\in\Ct^{n+1}$ is an isolated singular point of $s$. Let
$X_t=s^{-1}(t)\cap B_\epsilon$ for $0<|t|<\delta$ be the Milnor fiber
with the usual choice of $0<\delta\ll\epsilon\ll 1$.
For
a given germ at $0$ of a holomorphic $(n+1)$-form $\eta \in
\Omega^{n+1}_{\Ct^{n+1},0}$ define a quotient of forms by:
$$\left(\eta/ds\right)_{|X_t}=Res
\left(\frac\eta{s-t}\right)\in\Omega^n_{X_t}\,.$$
Let $\zeta_t\subset X_t$ be a continuous multivalued family of $n$-cycles
in the Milnor fibers. The function $$I^\eta_{\zeta}(t)=
\int_{\zeta_t}\eta /ds$$ is a holomorphic (multi--valued)
function. By \cite{Ma} or \cite{AGVII} \S13.1 the function
$I^\eta_{\zeta}(t)$ can be expanded in a series
$$I^\eta_{\zeta}(t)= \sum_{\alpha ,k} a_{\alpha ,k}t^\alpha(\log\,
t)^k\,,$$ where the numbers $\alpha$ are rationals greater then
$-1$ and $k$ are natural numbers or $0$.
When we consider all the possible families of cycles we
obtain so-called {\it geometric section} $S(\eta)$ of the cohomology
Milnor fiber.
We recall that cohomology Milnor fiber is a flat vector bundle equipped
with Gauss-Manin connection. Its fiber over $t$ is $H^n(X_t)$.
If we fix $t_0\neq 0$ we can write
$$S(\eta)=\sum_{\alpha ,k} A_{\alpha ,k}t^\alpha(\log\,t)^k,$$
with $A_{\alpha, k}\in H^n(X_{t_0})$.
The smallest
exponent $\alpha$ occurring in the expansion of $S(\eta)$ is
called the {\it order} of $\eta$. The smallest possible order among
all the forms $\eta$ is the order of $dz_0\wedge\dots\wedge dz_n$.

\begin{pr} \label{osi} Suppose that $X$ has isolated
singularities. Let $\omega\in\Omega^{n+1}_M( X )$ be a
meromorphic form with a first order pole along $X$.
If the order of $s\omega$ is greater than zero at each singular
point, then the residue class of $\omega$
lifts to intersection cohomology of $X$.\end{pr}

\begin{re}\rm For simple singularities with $n\geq 2$ the order of
any form is greater than zero.\end{re}

\Proof~~is based on the following easy local homological computation
(\cite{We2}, 2.1):

\begin{pr} \label{vli}If $X$ has isolated singularities then a differential
$n$-form on $X_\o$ defines an element in intersection cohomology if
and only if it vanishes in cohomology when restricted to the links of
the singular points.\end{pr}
Each cycle $\zeta_0$ in the link can be
extended to a family of cycles in the neighbouring fibers. We can
approximate the value of the integral $\int_{\zeta_0}\r$ by the
oscillating integral of $\eta=s\omega$. If all the exponents in
$I^\eta_{\zeta}(t)$ are greater than zero, then the limit integral
for $t=0$ vanishes. Therefore $[\r ]=0$ in the cohomology of each
link.
\qed

\begin{re}\rm Proposition \ref{osi} is a special case of the Theorem
\ref{lift}, although formulation of \ref{osi} is in terms of
oscillating integrals. By \cite{AGVII}, \S13.1 Th.1, the order of
$s\omega$ is greater than zero if and only if $\omega$ has canonical
singularities. Then $[\omega]\in W_{n+2}H^{n+1}(M\setminus X)$ and
\ref{lift} applies.\end{re}

\section{Quasihomogeneous isolated hypersurface singularities }
\label{qs}

More precise information about the exponents occurring in the
oscillating integrals can be obtained for isolated
quasihomogeneous singularities. All the simple singularities are
of this form. The resulting statement for the residue forms is
expressed in terms of weights. The weights of polynomials
considered here should not be confused with the weights in the
mixed Hodge theory. It is rather related to the Hodge filtration.
The relation is subtle and it will not be
discussed here. Let $a_0, a_1,\dots,a_n\in\Nt$ be the weights attached
to coordinates in which the function $s$ is quasihomogeneous. For
a meromorphic form of the top degree we compute the weight in the
following way:
$$v\left(\frac gsdz_0\wedge\dots\wedge
dz_n\right)=v(g)-v(s)+\sum_{i=0}^nv_i\,.$$

\begin{tm} \label{qua} Suppose that $X$ has isolated
singularities given by quasihomogeneous equations
in some coordinates. Let $\omega\in\Omega^{n+1}_M( X )$ be a
meromorphic form with a first order pole along $X$.
Suppose $\omega$ has no component of the weight 0 at each singular point.
Then the residue class of $\omega$
lifts to intersection cohomology of $X$.\end{tm}

\Proof. To apply Proposition \ref{vli} we will show that
$res(\omega)_{|L}=0$. It suffices to check that $\omega$ is
exact in
a neighbourhood of the singular points.
The calculation is local, we may assume that $M=\Ct^{n+1}$ and
$\omega\in\Omega^{n+1}_{\Ct^{n+1}}(X)$ is rational.
Suppose that $\omega$ is quasihomogeneous:
$$\omega={g\over s}dz_0\wedge\dots\wedge dz_n$$
with $g$ quasihomogeneous of degree $v(g)$. Then $g/s$ is
quasihomogeneous of degree $v(g)-v(s)$. This means that
$$\sum_{i=0}^n a_i{\partial(g/s)\over \partial z_i} z_i =
(v(g)-v(s)){g\over s}\,.$$
Let us define a form
$$\eta={g\over s}\sum_{i=0}^n (-1)^ia_i\,
z_i\,dz_0\wedge\dots^{ \vee\hskip -4pt^i}\hskip -2pt\dots\wedge
dz_n\,.$$
Then $$d\eta=\sum_{i=0}^na_i\left({\partial(g/s)\over \partial z_i}
z_i+{g\over s}\right)dz_0\wedge\dots\wedge dz_n=
\left(v(g)-v(s)+\sum_{i=0}^na_i\right){g\over s}dz_0\wedge\dots\wedge 
dz_n\,.$$
Therefore if $v(\omega)=v(g)-v(s)+\sum_{i=0}^na_i\neq 0$ then
$\omega ={1\over v(g)-v(s)+\sum_{i=0}^na_i}d\eta$.
\qed

\begin{re}\rm Conversely, if $\omega\neq 0$ is quasihomogeneous of degree 0
then the residue form restricted to the link $L$ of the singular
point is nonzero, $res(\omega)_{|L}\neq 0$. To see this consider the quotient 
$$L/S^1\subset \Pt(a_0,\dots a_n)$$
in the weighted projective space. 
Here $L$ is the link of the singular point; it is homeomorphic to the
intersection of $X$ with the unit sphere.
%$\{(z_0,z_1,\dots,z_n)\in\Ct^{n+1}:\sum|z_i|^{2a_i}=1\}$. 
The circle
acts on $\Ct^{n+1}$ diagonally with weights $a_i$.
Integrating along the fibers of the
quotient map one obtains a holomorphic form called by us {\it the second 
residue}
$$Res^{(2)}(\omega)=
\int_{S^1}Res(\omega)\neq 0\in \Omega^{n-1}_{L/S^1}\,.$$
Although $L/S^1$ is not smooth, it may have only quotient
singularities and the Hodge theory applies. Therefore
$[\int_{S^1}Res(\omega)]\neq 0\in H^{n-1}(L/S^1)$. We will
illustrate this construction by an example.
\end{re}

\begin{re}\rm Fix a real number $p>1$. If $v(\omega)>0$ one can construct
on $X_\o$ a conelike metric adapted to the quasihomogeneous
coordinates such that $|\omega|$ is integrable in the $p$-th
power. By \cite{We1} $L^p$-cohomology is isomorphic to
intersection cohomology for a perversity $\underline q$ with
${2n\over p}-1\leq \underline q(2n)< {2n\over p}$. For large $p$
it is isomorphic to cohomology of the normalization of $X$. This
way (again) we obtain an explicit lift to cohomology.\end{re}

\begin{ex}\rm Elliptic singularity: 
Consider a singularity of the type $P_8$ (elliptic singularity) in a form
$$s(z_0,z_1,z_2)=z_1^3+pz_0^2z_1+qz_0^3-z_0z_2^2$$
where $p$ and $q$ are real numbers
such that the polynomial $z^3+pz+q$ does not have double roots.
Let $$\omega={1\over s}dz_0\wedge
dz_1\wedge dz_2\,.$$
Then
$$\r=- {1\over 2 z_0z_2} dz_0\wedge dz_1$$
for $z_0z_2\neq 0$. The second residue is equal to
$$Res^{(2)}(\omega)={dz_1\over
z_2}={dz_1\over\sqrt{z_1^3+pz_1+q}}\,.$$ 
If we integrate $Res^{(2)}(\omega)$ along the real part of the
elliptic curve $L/S^1\subset\Pt^2$ we obtain the classical elliptic integral.
\end{ex}

\begin{re}\rm It would be enough to show that $Res^{(2)}(\omega)$ is 
nonzero as a form, since it is holomorphic it cannot vanish in
cohomology. Counting the homogenity degree it is immediate to check
that the second residue of the of the form ${1\over s}dz_0\wedge
dz_1\wedge dz_2$ is nontrivial for any homogeneous polynomial of the
degree 3. The coefficients of $s$ do not have to be real. Nevertheless we
find it interesting to see exactly what kind of 
numbers can appear as values of the second residue.
\end{re}
\vfil\sl

Andrzej Weber

Instytut Matematyki, Uniwersytet Warszawski

ul. Banacha 2, 02--097 Warszawa, Poland

e-mail: aweber@mimuw.edu.pl
\end{document}